**Note on Integer Factoring Methods III**
N. A. Carella, July 2007.

**Abstract:** The best unconditionally proven deterministic integer factorization algorithms have exponential running time complexities of $O(N^{1/4})$ arithmetic operations, and conditional on the Riemann hypothesis, there is a deterministic algorithm of exponential running time complexity $O(N^{1/5})$. This note proposes a new deterministic integer factorization algorithm of deterministic exponential time complexity $O(N^{1/6})$. Furthermore, an algorithm for decomposing composite integers that have factor differences of the form $q - p = (r - 1)N^{1/2} + u$, where $r \geq 1$ is a fixed parameter, and $|u| < N^{1/3+\varepsilon}$, in deterministic logarithmic time and various other results are included.

**1 Introduction**
The best unconditionally proven deterministic integer factorization algorithms have exponential running time complexities of $O(N^{1/4})$ arithmetic operations, for example, Coppersmith algorithm, the Strassen algorithm, the Pollard algorithm, and a few new ones, see [RL], [CP], [MP], [MZ]. Moreover, conditional on the Riemann hypothesis, there is the Shank algorithm of deterministic exponential running time complexity $O(N^{1/5+\varepsilon})$, $\varepsilon > 0$, [CP, p. 248]. This note proposes a new deterministic integer factorization algorithm of exponential time complexity $O(N^{1/6})$. Furthermore, an algorithm for decomposing composite integers that have factor differences of the form $q - p = (r - 1)N^{1/2} + u$, where $r \geq 1$ is a fixed parameter, and $|u| < N^{1/3+\varepsilon}$, in deterministic logarithmic time and various other results are included.

The content of this note is intended for a general audience and includes as much elementary details as possible. Exceedingly complex proofs are not included but references are provided. Section 2 recalls a few of the known works on divisors in residue classes and lattice reduction theory. The main contributions are Theorem 5 in Section 2, and Theorem 14 in Section 3. The last Section consists of (optional) background materials on lattice reduction theory.

**2 Residue Class Method**
The idea of expressing a composite integer as a product of linear equations dates back to 1800's or earlier. The algebraic equation

$$(mx + c)(my + d) = N, \tag{1}$$

where $m$ is a fixed modulo, was employed by early authors to factor and primality test contemporary large numbers, see [W1, p.102]. The analysis and application of (1) has been extended and improved by many authors, see [LH], [CG], [BL] etc.

A linear term $mx + c$ or equivalently $p \equiv c \mod m$ is called a *residue class*. A squarefree composite integer has a nontrivial representation as (1) with $x \neq y$.

For $m \approx N^{1/2}$, it is trivial to verify that an integer has at most 1 or 2 factors in a single residue class $mx + c$. But for $m \approx N^{1/4}$, it is more difficult. A careful inspection of squarefree integers reveals that there are at must 14 proper factors and at must 4 prime factors. In fact something similar holds for arbitrary integers.

**Theorem 1.** ([LH]) For every $\alpha \in \mathbb{R}$ with $\alpha > 1/4$ there exists a constant $c(\alpha)$ with the following properties. If $r, s, n$ are integers satisfying $n > 0$, $s > n^\alpha$, $\gcd(r, s) = 1$. Then the number of positive divisors of $n$ that are congruent to $r$ modulo $s$ is at most $c(\alpha)$.

On the other hand, for $m < N^{1/4}$ the number of factors of some integers $N$ that fall in a single residue class can be a slowly increasing function of $N$.

**Theorem 2.** ([LH]) Let $r, s$ and $n$ be integers satisfying $0 \leq r < s < n$, $s > n^{1/3}$, $\gcd(r, s) = 1$. Then there exists at most 11 positive divisors of $n$ that are congruent to $r$ modulo $s$, and there is a logarithmic time algorithm for determining all these divisors.

**Note 1:** The standard term *polynomial time* has been replaced with the more appropriate term *logarithmic time*. This is patterned after the closely related term *exponential time*.

The proof of this result consists of a series of equations similar to the Euclidean algorithm. A recent paper treats the concept of residue classes in term of lattice reduction theory, see [CG]. A related result was given earlier in [CR].

**Theorem 3.** ([CR]) If the $\log_2(N)/4$ least significant bits of the factor $p$ of $N = pq$, $p < q < 2p$, are known, then the factorization of $N$ has logarithmic time complexity.

Proof: Let $N$ be a balanced integer of $n = \log_2(N)$ bits. From the given the least significant $n/4$ bits part $x_0$ of $p = 2^{n/4}x + x_0$, use the congruence $x_0 y_0 \equiv N \mod 2^{n/4}$ to compute the least significant $n/4$ bits part $y_0$ of $q = 2^{n/4}y + y_0$, (similar calculation applies to the most significant parts). Then it follows that the polynomial equation

$$f(x, y) = (2^{n/4}x + x_0)(2^{n/4}y + y_0) - N = 2^{n/4}xy + y_0 x + x_0 y + (x_0 y_0 - N)/2^{n/4} = 0 \tag{2}$$

has a small root $(x_1, y_1)$ such that $0 \leq |x_1| \leq X < N^{1/4}$, $|y_1| \leq Y < N^{1/4}$. Moreover, the height

$$\|f(xX, yY)\|_\infty = \max\{2^{n/4}XY, y_0 X, x_0 Y, (x_0 y_0 - N)/2^{n/4}\} \approx N^{3/4} \tag{3}$$

of the associated polynomial $f(xX, yY)$ is sufficiently large, and by lattice reduction methods the roots $(x_1, y_1)$ such that $0 \leq |x_1 y_1| \leq XY < \|f(xX, yY)\|^{2/3} \approx N^{1/2}$ can be determined in deterministic logarithmic time, see Theorem 25. ∎

Given about 256 bits of a factor of a 1024-bit balanced integer, the current implementations of this algorithm can split the integer within minutes, see the literature. The new algorithm reduces it to about 170 bits, and the special case with $|x_0 y_0| < 2^{n/4+1}$ has an equivalent complexity of about 64 bits.

**Theorem 4.** Let $m \geq N^{1/4}$. If $cd < 2m$ then the factorization of the integer $N = (mx + c)(my + d)$ with two nearly equal factors has deterministic exponential time complexity $O(N^{1/16+\varepsilon})$, $\varepsilon > 0$.

Proof: Factor the integer $r \equiv N \bmod m$, and apply Theorem 3 to each pair $(c, d)$ such that $r = cd$. To estimate the time complexity, observe that the factorization of $r$ requires $O(m^{1/4}) = O(N^{1/16})$ arithmetic operations and there are at most $o(m^{\varepsilon})$ such pairs, $\varepsilon > 0$. Thus the entire process requires at most $O(N^{1/16+\varepsilon})$ arithmetic operations. ∎

Roughly speaking, the general principle involved in Theorem 4 shifts the factorization of the integer $N$ to the factorization of the smaller integer $cd < N^{2\alpha}$, $0 < \alpha < 1/2$. Theorems 2 and 3 work over the ranges $1/3 < \alpha < 1/2$, and $1/4 < \alpha < 1/2$ respectively. However, it probably can go as low as $1/6 < \alpha < 1/2$, see Theorem 8. A general heuristic algorithm for determining any residue classes $(c, d)$ is given in Subsection 2.1.

Since a brute force search can be used to determine either the residue class of $p$ or equivalently the least significant $\log_2(N)/4$ bits of $p$, Theorems 2 and 3 are general purpose integer factoring algorithms of deterministic time complexities $O(N^{1/3})$ and $O(N^{1/4})$ respectively. There are other algorithms of similar complexities, see [RL], [MP], [MZ], [OD], [BR], [S], [LA], [CP] and other sources.

In spite of many years of extensive research efforts to reduce the deterministic time complexity $O(N^{1/4})$, it remains untouched. Similarly, the conditional deterministic time complexity $O(N^{1/5})$ of the Shank algorithm has not been settled, see [CP, p. 248] for the current perspective.

**Theorem 5.** Let $N = pq$, $p < q < 2p$. If the $(1/6)\log_2(N)$ most significant bits of the factor $p$ are known, then the factorization of $N$ has logarithmic time complexity.

Proof: Let $P_0 = (1/6)\log_2(N)$ most significant bits of the factor $p$, and compute $Q_0 = (1/6)\log_2(N)$ most significant bits of the factor $q$. Put $Q_0 \pm a = Mz_0$ or a similar transformation, where $0 \leq |a| < 10$, and $1 \leq z_0 \leq O((\log N)^B)$, $B > 0$ constant. Rewrite the usual polynomial equation $(P_0 + x)(Q_0 + y) - N = 0$ as

$$f(x, y, z) = (P_0 + x)(Mz + y) - N = xy + Mxz + P_0 y + MP_0 z - N = 0. \qquad (4)$$

Clearly, this is an irreducible polynomial over the integers $\mathbb{Z}$, and it has a small solution $(x, y, z) = (x_0, y_0, z_0)$, where $0 \leq |x|, |y| < N^{1/3+\varepsilon}$. Now apply a lattice reduction technique to construct another polynomial $g(x, y, z)$ of sufficiently small norm $\| g(xX, yY, zZ) \|_2$, which satisfies the inequalities

$$0 < \| g(xX, yY, zZ) \|_2 < \| f(xX, yY, zZ) \|_\infty = N. \qquad (5)$$

By Corollary 27, $g(x, y, z)$ is not a multiple of $f(x, y, z)$, so these polynomials are algebraically independent over the polynomials ring $\mathbb{Z}[x,y,z]$. Moreover, since $f(x, y, z)$ and $g(x, y, z)$ share the same small root $(x, y, z) = (x_0, y_0, z_0)$ and $z_0$ is known, the solution $(x, y) = (x_0, y_0)$ is recovered by means of the resultant

$$\text{Res}_x(f(x, y, z_0), g(x, y, z_0), y) \quad \text{or} \quad \text{Res}_y(f(x, y, z_0), g(x, y, z_0), x). \qquad (6)$$

The complexity of the algorithm is deterministic logarithmic time as the lattice reduction technique. ∎

Replacing a constant with a known transformation offers several advantages: (i) it circumvents the requirement for a third algebraically independent polynomial, (ii) it expands the ranges of the variables $x$ and $y$, (iii) it works with either the most significant bits or the least significant bits. Depending on the lattice construction, a significantly better result is quite possible. The technical details for constructing polynomial lattices appear in [EJ] and similar sources, the latest attempt to construct third algebraically independent polynomials appears in [BA].

### 2.1 Determination of Residue Classes

Very little work has been done on the calculations of the residue classes per se. The well known result below uses an interactive brute force search to compute the residue classes, it interacts with other algorithms to determine the correct residue classes $p \equiv c \bmod m$, $q \equiv d \bmod m$ of $N$. Thus it is not very effective whenever $m = N^\alpha$, $\alpha > 0$.

***Proposition 6.*** Let $m = O(N^\alpha)$ be a fixed modulo. Then the constants $c$ and $d$ are computable in $O(N^\alpha)$ operations modulo $N$.

Proof: For each $0 \neq c \in \mathbb{Z}_m$ such that $\gcd(c, m) = 1$, there exists a unique $d$ such that $cd \equiv N \bmod m$, which is computable via the Euclidean algorithm, and multiplications modulo $N$. ∎

The early details of a different technique for computing the residue classes is introduced here. Although the current version uses a brute force search, the resulting algorithm can be effective because it does not interact with other algorithms. It is a one-time-calculation stand-alone-algorithm. The algorithm accepts an integer $N$ and a prime number $m$, and it generates a set of residue class pairs $R = \{ (c, d) : cd \equiv N \bmod m \}$ that contains the residue classes of $N$ modulo $m$.

*Algorithm* I
Input: Integer $N$, and prime number $m$.
Output: A set $R = \{ (c, d) : cd \equiv N \bmod m \}$ of probable residue classes.
1. For each $x \geq 0$, or in the range $2\sqrt{r_0} \leq x < 2m$, solve $x^2 - 4cd \equiv y^2 \bmod m$, where $cd = r_0 + r_1 m$, and $r_0 \equiv N \bmod m$.

2. Put $c = (x + y)/2$ and $d = (x - y)/2$, then collect the pairs $(c, d)$ that satisfy the sieving relations $cd \equiv N \bmod m$, and $cd \equiv cd \bmod m^2$ in the set $R$.
3. Return $R$.

The $m$ prime restriction simplify the algorithm and eliminates difficult cases. In addition, since $0 \leq c, d < m$, and the sum $0 < c + d < 2m$, all the integer solutions $(x, y) = (c + d, c - d)$ of the equation $x^2 - 4cd = y^2$ can be computed via the congruence equation $x^2 - 4cd \equiv y^2 \bmod m$, and a lift $x^2 - 4cd \equiv y^2 \bmod m^2$. Only those integer solutions $(x, y) = (c + d, c - d)$ that satisfy the sieving relations $cd \equiv N \bmod m$, and $cd \equiv cd \bmod m^2$ are collected in the set of possible residue classes $R$.

Any of the advanced techniques for solving the equation $x^2 - 4cd \equiv y^2 \bmod m^v$, $v \geq 1$, are applicable here. Some special cases possible are as follows:
(1) If $cd = a^2$, then the congruence $x^2 - 4cd \equiv y^2 \bmod m^2$ has a unique solution.
(2) If $cd < m$, then the solution of the congruence $x^2 - 4cd \equiv y^2 \bmod m$ holds over the integers.
The cardinality of the set $R$ is $\#R = O(\log N)$ on average and it is $\#R = o(N^\varepsilon)$ on a worst case condition, $\varepsilon > 0$. These estimates are derived from the average order and the asymptotic behavior

$$\sum_{n \leq x} v(n) = x \log N + E(x) \quad \text{and} \quad \lim_{n \to \infty} \frac{v(n)}{n^\varepsilon} = 0 \tag{7}$$

respectively of the divisor counting function $v(n) = \#\{ d : d \mid n \}$, the error term $E(x) < 5x$.

## 3 Quasi Residue Class Method
This section starts with several known but interesting results, and a probably new perspective in integer factorization.

**3.1 Landry-Pepin method.** The sketch provided is a minor generalization of the ideas exploited by Landry, Pepin, Lehmer, et cetera. The early success of this technique in the factorization of $2^{64} + 1$ is described in [W2]. The reader should consult [B, lxiv], [CP], et cetera for other related results in both primality testing and integer factorization.

Let $m$ and $n \leq N^\alpha \leq N^{1/2}$ be fixed moduli. The factors of a composite integer $N$ are expressed as noncongruent residue classes

$$p = mx + c \qquad \text{and} \qquad q = ny + d, \tag{8}$$

where the variables $x$ and $y$ satisfy $0 \leq |x|, |y| < 2N^{1/2-\alpha}$, and the integers $c$ and $d$ satisfy $0 \leq |c|, |d| < N^\alpha$. The corresponding algebraic equation has the form

$$mnxy + dmx + cny + cd = N. \tag{9}$$

There are several approaches to proceed from equations (8) and (9) to other equations that are relevant to integer factorization algorithms. Some of the details of these approaches are sketched below. Rewrite the last expression as

$$mnxy + dp + cq = N + cd. \tag{10}$$

**Proposition 7.** Let $m, n = N^\alpha \leq N^{1/4}$ be fixed moduli, and suppose that $0 \leq |c|, |d| \leq N^\beta < N^{1/4}$. Then the factors of the integer $N = pq$, $p < q < 2p$, are computable in $O(N^{1/2-2\alpha+3\beta})$ arithmetic operations.

Proof: Equation (10) implies that the scaled sum of primes $z_t = dp + cq$ is also given by

$$z_t = dp + cq = mnt + z_0, \tag{11}$$

where

$$z_0 \equiv N + cd \mod mn \quad \text{and} \quad t \in \mathbb{Z}, \tag{12}$$

(here the sign of $z_0$ is the same as the sign of $dp + cq$). The factors are recovered as rational roots $\pm p$ and $\pm cq/d$ of the polynomial $dX^2 \pm z_t X \pm cN = 0$.

To estimate the computational time complexity, observe that the penultimate identity implies that

$$|mnt + z_0| = |dp + cq| \leq (p+q)N^\beta < 3N^{1/2+\beta}. \tag{13}$$

Therefore the parameter $t$ satisfy $0 \leq |t| < N^{1/2-2\alpha+\beta}$. The extra term $N^{2\beta}$ in the stated time complexity estimate accounts for the cost of determining the correct pair $c, d$ modulo $m$ and modulo $n$ respectively, see Proposition 6. ∎

The special case of $m = n$ and $c = d = 1$, due to Lehmer, see [BL], has a deterministic time complexity of $O(N^{1/2-2\alpha})$ arithmetic operations. More generally, if $p = mx + c$ and $q = my + d$, with $m \geq N^{1/4}$ and $0 \leq |c|, |d| \leq O((\log N)^A)$, then the integer $N = pq$ can be factored in deterministic logarithmic time $O((\log N)^A)$, $A > 0$ constant.

Fix a modulo $m \geq N^{1/4}$. If there exists a pair $c, d$, then a brute force search can be performed to determine a partial factorization $N - cd = mR$. But the density of such integers is negligible. On the other hand, a search of the parameters space $\{c\} \times \{d\}$ for a suitable modulo $m \geq N^{1/4}$ is a significantly more effective strategy for factoring an integer since the density of a variable modulo $m = m(c, d)$ is considerably larger.

Other special cases are readily derived from Proposition 7, in particular, the use of a pair of relatively prime moduli reduces the moduli sizes from $N^{1/4}$ to $N^{1/6}$.

**Theorem 8.** Let $m$ and $n \geq N^{1/6}$ be relatively prime moduli and let $0 \leq |c|, |d| \leq O((\log N)^A)$. Then the composite integer $N = pq$, $p < q < 2p$, can be factored in deterministic logarithmic time $O((\log N)^{2A})$, $A > 0$ constant.

Note on Integer Factoring Methods III

The proof uses the fact that gcd($m, n$) = 1 to construct an irreducible polynomial $f(x, y) = (mx + c)(ny + d) - N$ of sufficiently large height $\| f(xX, yY) \| \approx N$, and lattice reduction methods to determine the small roots of the equation $f(x, y) = 0$, as in Theorem 3. This is quite effective, however, very few integers satisfy these conditions. In fact, the density of these integers is approximately

$$\frac{(\log N)^A}{\varphi(m)} \times \frac{(\log N)^A}{\varphi(n)} \approx \frac{(\log N)^{2A}}{N^{2/3}}, \tag{14}$$

which is essentially zero as $N$ increases. The density estimate is derived from the Prime Number Theorem, and Dirichlet Theorem on primes in arithmetic progression $\{ ax + b : \gcd(a, b) = 1, x \geq 0 \}$, these are

$$\pi(x) = \frac{c_0 x}{\ln x} + O(xe^{-c_1\sqrt{\ln x}}) \quad \text{and} \quad \pi(x, a, b) = \frac{c_0 x}{\varphi(a)\ln x} + O(xe^{-c_1\sqrt{\ln x}}) \tag{15}$$

respectively, see [S, p86], [E, p. 266] and similar sources for discussions. It is also assumed that the factors $p$ and $q$ are independent random variables in the ranges $\sqrt{N/2} \leq p \leq \sqrt{N}$ and $\sqrt{N} \leq q \leq \sqrt{2N}$ respectively.

**Proposition 9.** Let $m = O((\log N)^A)$, $A > 0$ constant, and let $\varepsilon > 0$. Then the linear Diophantine equation $dx + cy - mz = b$ has a small solution in the ranges $N^{1/2-\varepsilon} < |x|, |y| < N^{1/2}$ and $0 \leq |z| < N^{1/2}$.

Here the coefficients $b$, $c$, and $d$ are known. To prove this claim, consider the $m$-expansion of the integer

$$N = a_2 m^2 + a_1 m + a_0 = am^2 + bm + cd. \tag{16}$$

The coefficients of the expansion are matched to the corresponding terms in the algebraic equation

$$m^2 xy + (dx + cy)m + cd = N. \tag{17}$$

However, due to the nonuniqueness of the expansion for $m < N^{1/3}$, this relationship requires a third variable $z$, which accounts for the carry in the expansion. Specifically

$$xy = a - z \quad \text{and} \quad dx + cy = b + mz. \tag{18}$$

Proper selection of the modulo as $m = O((\log N)^A)$ leads to small solution in a 4-dimensional integers lattice $\mathbb{Z} \times \mathbb{Z} \times \mathbb{Z} \times \mathbb{Z} \times \mathbb{Z}$. In light of the multidimensional Euclidean algorithm, integer relations methods, and the like, see [FF], [FH], [AK], [HL], etc., this is of interest in integer factoring algorithms.

**3.2 Fermat Difference of Squares Method.** The difference of squares method factors an odd integer $N = pq$ in a few arithmetic steps whenever there is a pair of factors $p$ and $q$ close to the geometric mean $\sqrt{N} = \sqrt{pq}$, or equivalently whenever the factors difference $q - p \leq O(N^{1/4})$ is small, [RL, p. 147].

**Proposition 10.** (1) An integer $N$ is represented as $N = x^2 - y^2$ if and only if $N \neq 4M + 2$.
(2) A prime $N > 2$ has a unique representation as a difference of consecutive squares.
(3) The number of solutions is $O(\log N)$ on average and $o(N^\varepsilon)$ solutions asymptotically.

Proof: In the first statement use the fact that $x^2 - y^2 \equiv 0, 1, 3 \mod 4$, and in the second statement use $N = (2^{e-1}m + 1)^2 - (2^{e-1}m)^2$ or $(2^{e-1}m)^2 - (2^{e-1}m - 1)^2$, where $N = 2^e m \pm 1$, $m$ odd, and $e \geq 1$. The last statement follows from (7). ∎

Observe that any integer $N \geq 1$ has a representation as

$$4N = x^2 - y^2, \qquad (19)$$

and that any solution of (19) is of the form $x = p + q$, $y = q - p$ with $p, q \mid N$. The extreme solution $x = N + 1$, $y = N - 1$ does not lead to a nontrivial factorization of $N$, so it is viewed as the trivial solution. The factors of an arbitrary integer $N$, which can be composites or primes, vary from $p = q = N^{1/2}$, to $p = N/2$, $q = 2$, and $p = N$, $q = 1$. Equality occurs if and only if $N$ is a square.

A prime number $N > 2$ has a unique solution of large consecutive integers, which is the trivial solution. But if $N$ is not prime, then there is a nontrivial solution such that $x \geq 2N^{1/2}$ is an integer in the sequence of integers

$$x_0 = 2\sqrt{N}, \; x_1 = 2\sqrt{N} + 1, \; x_2 = 2\sqrt{N} + 2, \; ..., \; x_n = (N+4)/2. \qquad (20)$$

Technically $x_i = [2\sqrt{N}] + i$, where the bracket $[x]$ is the largest integer function, however the bracket is often omitted to simplify the notation.

**Proposition 11.** Let $\varepsilon > 0$. If $p \mid N$ and $|p - N^{1/2}| < N^{1/4+\varepsilon}$ then the difference of squares method decomposes the integer $N$ in deterministic logarithm time.

Proof: Put $p + q = p + N/p = x$. Starting at $x = 2N^{1/2}$ as in the sequence (20), it easy to show that the number of steps required to find a solution is

$$p + N/p - 2\sqrt{N} = (p - \sqrt{N})^2 / p. \qquad (21)$$

Clearly, if there is a factor $p$ sufficiently close to $N^{1/2}$ the procedure is successful, and runs in deterministic logarithmic time (even constant time). Otherwise,

$$p \leq N^{1/2} - N^{1/4+\varepsilon} \quad \text{and} \quad N^{1/2} + N^{1/4+\varepsilon} \leq q \qquad (22)$$

and the algorithm runs in exponential time $O(N^{2\varepsilon})$. ∎

The previous analysis does not utilize any special properties of the factors of $N$ beside the distance $|q - p|$ between the factors. However, in some cases there are acceleration techniques that utilize other properties of the integer $N$ such as exclusion moduli, etc, to speed up the process. An interesting acceleration technique is effected by the triangular number identity

$$1^3 + 2^3 + 3^3 + \cdots + k^3 = (k(k+1)/2)^2. \qquad (23)$$

This identity facilitates a search for a solution $(x, y)$ of $4N = x^2 - y^2$ at every $2N^{1/4}$ other number in the sequence of square numbers

$$\begin{aligned}
x_0^2 &= (m(m+1)/2)^2, \qquad m = [2N^{1/4}], \\
x_1^2 &= x_0^2 + (m+1)^3, \\
x_2^2 &= x_1^2 + (m+2)^3, \\
&\ldots \\
x_i^2 &= x_{i-1}^2 + (m+i)^3,
\end{aligned} \qquad (24)$$

instead of a consecutive search as in the standard algorithm (20). This acceleration technique works whenever the sum $p + q = k(k + 1)/2$ (or the difference $q - p$) is a triangular number. The estimated starting point is taken from the inequality $k(k + 1)/2 = x \geq 2N^{1/2}$. For example, to factor $N = 23 \cdot 113$ this technique uses only 3 steps:

$$\begin{aligned}
x_0^2 &= (m(m+1)/2)^2 = 105^2, \quad x_0^2 - 4N \neq \text{square}, \quad m = [2N^{1/4}] = 14, \\
x_1^2 &= x_0^2 + (m+1)^3 = 120^2, \quad x_1^2 - 4N \neq \text{square} \\
x_2^2 &= x_1^2 + (m+2)^3 = 136^2, \quad x_2^2 - 4N = y_2^2 = 90^2,
\end{aligned} \qquad (25)$$

as opposed to 35 steps in the standard algorithm.

It is quite possible that this acceleration technique can also be used to factor integers that have nontriangular sums $p + q \neq k(k + 1)/2$. In this case the sequence of squares (24) is used as an approximation of $(p + q)^2$.

**3.4 Uniform Difference of Square Method.** The difference of square method excels in the factorization of integers $N = pq$ composed of factors of the form $p < N^{1/2}$ and $q = p + u$, with $|u| < N^{1/4}$. A new extension of this technique achieves the same performance for any ratio $r = q/p \geq 1$ uniformly. A proof almost identical to the case $r = 1$ as in Proposition 11 is workable. However, one based on the more instructive and general concept of lattice reduction methods is given.

***Theorem* 12.** Let $N = pq$ and suppose that $p < N^{1/2}$ and $q = rp + u$, with $|u| < N^{1/4}$. If the ratio $r = q/p \geq 1$ is given, then the factorization of the integer $N$ has deterministic logarithmic time complexity.

Proof: Use the equation $N = p(rp + u)$ and the power series expansion of the function $f(x) = \sqrt{1+x}$ to derive suitable expressions for the factors $p$ and $q$. Specifically write

$$\sqrt{N} = p\sqrt{r}\sqrt{1+u/rp} = p\sqrt{r}\left(1 + \frac{u}{2rp} - \frac{u^2}{8r^2p^2} + \frac{3u^3}{16r^3p^3} - \cdots\right). \tag{26}$$

Simple algebraic manipulations yield

$$p = \sqrt{N/r} - \left(1 + \frac{u}{2r} - \frac{u^2}{8r^2p} + \frac{3u^3}{16r^3p^2} - \cdots\right) = \sqrt{N/r} + x, \tag{27}$$

$$q = \sqrt{rN} - \left(1 + \frac{u}{2} - \frac{u^2}{8rp} + \frac{3u^3}{16r^2p^2} - \cdots\right) = \sqrt{rN} + y,$$

where $0 \leq |x|, |y| < N^{1/4}$. Let $X, Y = N^{1/4}$, and let

$$f(x,y) = \left(\sqrt{N/r} + x\right)\left(\sqrt{rN} + y\right) - N = 0. \tag{28}$$

Then the height of the associated polynomial $f(xX, yY)$ satisfies the relation

$$\|f(xX, yY)\| = \max\left\{|XY|, |\sqrt{rN}X|, |\sqrt{N/r}Y|, |c_0|\right\} = \sqrt{rN}X \approx N^{3/4}, \tag{29}$$

where $c_0 = [\sqrt{N/r}][\sqrt{rN}] - N$. Therefore, the small solutions $0 \leq |x|, |y| < N^{1/4}$ of the polynomial $f(x, y)$ can be determined in deterministic logarithmic time, see Theorem 23. ∎

This is an effective procedure for factoring any integer $N = pq$ that satisfies the constraints

$$\left|p - \sqrt{N/r}\right| < N^{1/4+\varepsilon}, \left|q - \sqrt{rN}\right| < N^{1/4+\varepsilon} \tag{30}$$

whenever $r \geq 1$ is known. The case $r = 1$ reduces to a lattice reduction theory formulation of the difference of square method for the integers $N = pq$ such that $\left|p - \sqrt{N}\right| < N^{1/4+\varepsilon}, \left|q - \sqrt{N}\right| < N^{1/4+\varepsilon}$. The corresponding equation is $\left([N^{1/2}] + x\right)\left([N^{1/2}] + y\right) - N = 0$, where $0 \leq |x|, |y| < N^{1/4}$. The case $r = 2$ (or any other obvious choice of $r \geq 1$) is a new effective technique for factoring the integers $N = pq$ such that $\left|p - \sqrt{N/2}\right| < N^{1/4+\varepsilon}, \left|q - \sqrt{2N}\right| < N^{1/4+\varepsilon}$. The corresponding equation is $\left(\sqrt{N/2} + x\right)\left(\sqrt{2N} + y\right) - N = 0$.

An improved version of this new technique exploits the idea introduced in Theorems 5 and 13 to extend the ranges

$$|p - P_0| < N^{1/3+\varepsilon}, |q - Q_0| < N^{1/3+\varepsilon}, \tag{31}$$

where $P_0 = \sqrt{N/r}$ and $Q_0 = \sqrt{rN}$. The closest and best result related to this is the one described in [CG, Theorem 2.1]:

**Theorem 13.** ([CG])  Given $m$ and $n$ with $m = n^\alpha$, all $x$ such that $(m + x)$ divides $n$ and $|x| < n^\gamma$ can be found in polynomial time whenever

$$\gamma h(h-1) - 2u\alpha h + u(u+1) \leq -\varepsilon < 0,$$

for some integers $h > u > 0$ and some $\varepsilon > 0$. The largest value of $\gamma$ for which this can hold is $\alpha^2 - \varepsilon$.

The new result stated here has a wider range by a factor of $N^{1/12}$.

**Theorem 14.**  Let $N = pq$ and suppose that $p < N^{1/2}$ and $q = rp + u$, with $|u| < N^{1/3}$. If the ratio $r = q/p \geq 1$ is given, then the factorization of the integer $N$ has deterministic logarithmic time complexity.

Proof: Let $P_0 = \sqrt{N/r}$ and $Q_0 \pm a = \sqrt{rN} \pm a = Mz_0$ with $0 < |a|$, $z_0 \leq O((\log N)^B)$, $B > 0$ constant. Then there is an irreducible polynomial

$$f(x, y, z) = (P_0 + x)(Mz + y) - N = xy + Mxz + P_0 y + MP_0 z - N = 0, \tag{32}$$

over the set of integers $\mathbb{Z}$. Now apply a lattice reduction technique to obtain another polynomial $g(x, y, z)$ with the same small root $(x, y, z) = (x_0, y_0, z_0)$, where $0 \leq |x_0|, |y_0| < N^{1/3+\varepsilon}$. Then the root $(x, y) = (x_0, y_0)$ is recovered by means of the resultant

$$R_x(f(x, y, z_0), g(x, y, z_0), y) \quad \text{or} \quad R_y(f(x, y, z_0), g(x, y, z_0), x). \tag{33}$$

The complexity of the algorithm is deterministic logarithmic time as the lattice reduction technique. ∎

This result is equivalent to Theorem 5. In fact, an approximation of the ratio $r = q/p$ of just $\log(N)/6$ decimal places is sufficient.

**Selection of the Ratio $r$.** An inspection of the relations $N = pq$ and $q/2 \leq p \leq q \leq 2p$, immediately gives the inequalities

(i) $\sqrt{N/2} \leq p \leq \sqrt{N} \leq q \leq \sqrt{2N}$, (34)
(ii) $2\sqrt{N} \leq p + q \leq (3\sqrt{2}/2)\sqrt{N}$,
(iii) $0 \leq q - p \leq (\sqrt{2}/2)\sqrt{N}$.

These inequalities serve as a guide in the selection of the ratios $(r_0, s_0)$. Specifically,

(i') $\sqrt{2}/2 \le r_0 \le 1$, and $1 \le s_0 \le \sqrt{2}$ (35)

(ii') $2 \le r_0 + s_0 \le 3\sqrt{2}/2$,

(iii') $0 \le s_0 - r_0 \le \sqrt{2}/2$.

Similar inequalities as in (34) and (35) hold for the more general relations $N = pq$ and $aq \le p \le q \le bp$, where $a, b > 0$.

A large scale computation could involves a large grid of ordered pairs $(r_i, s_i)$, $0 \le i \le K$. As an example, a 21-point grid of ordered ratios $(r_i, s_i)$ is listed in the table below. These were obtained by uniform subdivisions of the interval $[.707, 1]$ into 21. Each pairs $(r_i, s_i)$ satisfies the constraint $r_i s_i = 1$.

| $r_i$ | $s_i$ | $r_i$ | $s_i$ | $r_i$ | $s_i$ |
|---|---|---|---|---|---|
| 0.707 | 1.414427 | 0.804667 | 1.242751 | 0.902333 | 1.108238 |
| 0.720952 | 1.387054 | 0.818619 | 1.221569 | 0.916286 | 1.091363 |
| 0.734905 | 1.360721 | 0.832571 | 1.201098 | 0.930238 | 1.074994 |
| 0.748857 | 1.335368 | 0.846524 | 1.181302 | 0.94419 | 1.059108 |
| 0.76281 | 1.310943 | 0.860476 | 1.162147 | 0.958143 | 1.043686 |
| 0.776762 | 1.287396 | 0.874429 | 1.143604 | 0.972095 | 1.028706 |
| 0.790714 | 1.264679 | 0.888381 | 1.125643 | 0.986048 | 1.01415 |

## 4 Theoretical Foundations

This optional section presents a limited introduction to lattice reduction theory and polynomial equations.

### 4.1 Lattice Reduction Theory

A subset of vectors $\{ v_1, v_2, \ldots, v_n \} \subset \mathbb{R}^n$ is a basis of a lattice $L$ if the linear equation

$$x_1 v_1 + x_2 v_2 + \cdots + x_n v_n = 0 \qquad (36)$$

has a unique solution $(0,\ldots,0)$. Etly, the subset of vectors $L = \{ x_1 v_1 + x_2 v_2 + \cdots + x_n v_n : x_i \in \mathbb{Z} \}$ in $n$-dimensional space $\mathbb{R}^n$ spanned by a basis $\{ v_1, v_2, \ldots, v_n \}$ is called a lattice $L$. Every lattice $L$ is a discrete subgroup of $\mathbb{R}^n$.

Any two basis $\{ u_1, u_2, \ldots, u_n \}$ and $\{ v_1, v_2, \ldots, v_n \}$ of the lattice $L$ are equivalent up to a linear change of variables $v_i = a_{i,1} u_1 + a_{i,2} u_2 + \cdots + a_{i,n} u_n$, where $A = ( a_{i,j} ) \in SL_n(\mathbb{Z})$ is an $n \times n$ integer matrix. The unimodular integer matrix $A = ( a_{i,j} )$ has unit determinant $\det(A) = \pm 1$.

The *discriminant* Disc($L$) of the lattice $L$ is defined by the determinant $\det(L) = \det[v_1, v_2, \cdots, v_n]$ of the corresponding $n \times n$ matrix $[v_1, v_2, \cdots, v_n]$. The determinant coincides with the volume of the *fundamental domain* $\mathcal{F} = \{\, x_1 v_1 + x_2 v_2 + \cdots + x_n v_n : 0 \leq x_i \leq 1 \,\}$ of the lattice.

***Theorem* 15.** [Hadamard 1949?] If $v_1, v_2, \ldots, v_n$ is a basis of lattice $L$, then $\mathrm{Disc}(L) \leq \|v_1\| \cdot \|v_2\| \cdots \|v_n\|$.

***Theorem* 16.** [Hermite 1879?] Let $L$ be a lattice of dimension $n \geq 1$. Then
(i) $L$ contains a nonzero vector $v$ of norm $\|v\| \leq \gamma_n \mathrm{Disc}(L)^{1/n}$, where $\gamma_n$ is a lattice invariant.
(ii) $L$ has a basis such that $\|v_1\| \cdot \|v_2\| \cdots \|v_n\| \leq \gamma_n^{n/2} \mathrm{Disc}(L)$.

The Hermite constant is the expression $\gamma_n = \left(\dfrac{\delta_n}{V_n}\right)^{2/n}$, where

$$V_n(S) = \int_{x_1^2 + \cdots + x_n^2 \leq r} dx_1 dx_2 \cdots dx_n = \frac{\pi^{n/2} r^n}{\Gamma(n/2 + 1)} = \begin{cases} \dfrac{\pi^m}{m!} & \text{if } n = 2m, \\[2mm] \dfrac{2^n \pi^m m!}{n!} & \text{if } n = 2m+1, \end{cases} \quad (37)$$

is the volume of a sphere $S$ of radius $r > 0$ in $n$-dimensional space, and $\delta_n$ is the maximal density of unit sphere packing.

***Theorem* 17.** (Blichfeldt 1914) The maximal density $\delta_n$ of unit sphere packing in $n$-dimension space $\mathbb{R}^n$ (lattice) satisfies $\delta_n \leq \dfrac{n+2}{2^{(n+2)/2}}$.

***Theorem* 18.** For $n > 8$, the constant $\gamma_n$ satisfies the asymptotic behavior

$$\frac{1}{2\pi e} \leq \frac{\gamma_n}{n} \leq \frac{1.7}{2\pi e}, \quad (38)$$

where $e = 2.7182818284\ldots$ is the logarithmus naturalis base.

The exact values for the first 8 cases are known: $\gamma_1 = 1$, $\gamma_2 = 4/3$, $\gamma_3 = 2$, $\gamma_4 = 4$, $\gamma_5 = 8$, $\gamma_6 = 64/3$, $\gamma_7 = 64$ and $\gamma_8 = 256$. The sphere packing problem has a half a millennium tradition. The proofs for the known cases $\delta_1, \delta_2, \delta_3, \delta_4, \delta_5, \delta_6, \delta_7,$ and $\delta_8$ are the products of half a millennium of research, [CS].

The $i$th successive minimum $\lambda_i(L)$ of a lattice is the smallest real number for which there are $i$ linearly independent vectors $v_1, v_2, \ldots, v_i$ such that $\|v_j\| \leq \lambda_i(L)$ for $j = 1, 2, \ldots, i$. Specifically, $\|v_1\| \leq \lambda_1(L)$, $\|v_2\| \leq \lambda_2(L)$, $\ldots$ .

***Gaussian Estimates* 19.** The number of lattice points in a set $S \subset L$ is approximately

(i) $\dfrac{Vol(S)}{Vol(L)}$, (39)

(ii) $\lambda_1 = \min_{v \in L}\{\|v\|\} \leq \sqrt{\dfrac{n}{2e\pi}} Disc(L)^{1/n}$,

where $Vol(X)$ is the volume of the set $X \subseteq L$.

The volume of a lattice is the same as the volume of a fundamental domain $\mathcal{F}$, so a subset $S$ that contains $\mathcal{F}$ contains a vector $v \in L$.

**The Nearest Vector Problem.** Given a vector $v \in \mathbb{R}^n$, the nearest vector problem asks for the vector $w \in L$ that satisfies the relation

$$\min_{w \in L}\{\|v-w\|: v \notin L\},$$ (40)

where the distance function $\|v-w\|^2 = (a_1 - b_1)^2 + (a_2 - b_2)^2 + \cdots + (a_n - b_n)^2$, and the vectors are $v = (a_1, a_2, \ldots, a_n)$, and $w = (b_1, b_2, \ldots, b_n)$ respectively.

**The Shortest Vector Problem.** The shortest vector problem asks for the vector $v \in L$ of smallest norm. Specifically it satisfies the relation

$$\min_{v, w \in L}\{\|v-w\|: v \neq w\}.$$ (41)

**Orthogonal Bases And The Gramm Schmidt Process**
The Gramm-Schmidt process is a general purpose orthogonalization procedure in any inner product space. It converts an arbitrary basis of an inner product space into an orthogonal basis. Let $\{v_1, v_2, \ldots, v_n\}$ be a basis of integer vector of a lattice $L$. The G-S process generates the orthogonal basis $\{v_1^*, v_2^*, \ldots, v_n^*\}$ of rational vectors. The orthogonal basis is constructed as series of adjustments to the given basis as follows:

$$v_1^* = v_1,$$ (42)

$$v_2^* = v_1 - \dfrac{<v_1, v_1^*>}{\|v_1^*\|^2} v_1^*,$$

$$v_i^* = v_i - \sum_{j=1}^{i-1} \dfrac{<v_1, v_j^*>}{\|v_j^*\|^2} v_j^* = v_i - \sum_{j=1}^{i-1} \mu_{i,j} v_j^*,$$

for $i = 1, 2, \ldots, n$. Here the conversion coefficients $\mu_{i,j} = \dfrac{<v_i, v_j^*>}{\|v_j^*\|^2}$ for $\|v_j^*\| \neq 0$, else $\mu_{i,j} = 0$.

The term $\dfrac{<v, v_i>}{<v_i, v_i>} v_i = \dfrac{<v, v_i>}{\|v_i\|^2} v_i$ is the orthogonal projection of the vector $v$ onto the vector $v_i$.

The normalization $e_i = \dfrac{v_i^*}{<v_i^*, v_i^*>}$ spawns an orthonormal basis $\{ e_1, e_2, \ldots, e_n \}$ of the lattice. The important delta function relation $<e_i, e_j> = \delta_{i,j}$ characterizes an orthonormal basis.

Given an orthogonal basis $\{ v_1, v_2, \ldots, v_n \}$ and a vector $v \in L$, the vector has the orthonormal expansion

$$v = \frac{<v, v_1>}{\|v_1\|^2} v_1 + \frac{<v, v_2>}{\|v_2\|^2} v_2 + \cdots + \frac{<v, v_n>}{\|v_n\|^2} v_n. \tag{43}$$

***Theorem* 20.** ([L]) A basis $b_1, b_2, \ldots, b_w$ of a lattice $L$ is called a reduced basis if it satisfies the inequalities

$$\|b_1\| \leq \|b_2\| \leq \cdots \leq \|b_i\| \leq 2^{w(w-1)/4(w-1-i)} \det(L)^{1/(w-1-i)}, \tag{44}$$

for $i = 1, 2, \ldots, w$. The basis reduction algorithm uses $O(w^4 \log(B))$ arithmetic operations on integers of binary length $O(w\log(B))$, where $\|b_i\| \leq B$, $B \in \mathbb{R}$ is a constant and $\|\ \|$ is the standard norm.

The most recent algorithmic developments in lattice reduction methods are not considered here, but these are recommended to implement efficient algorithms.

### 4.2 Polynomials

Let $\mathbf{R}$ be a ring and let $\mathbf{R}[x_1 \ldots x_n] = \{ f(x_1 \ldots x_n) = \sum_{e=(e_1,\ldots,e_n)} a_e x_1^{e_1} x_2^{e_2} \cdots x_n^{e_n} : e_i \geq 0, a_e \in \mathbf{R} \}$ be the ring of polynomial functions of $n$ variables. For the application envisioned here the polynomial will be written in term of a single variable.

Let $f(x_1,\ldots,x_n) = a_k x_1^k + a_{k-1} x_1^{k-1} + \cdots + a_1 x_1 + a_0$, $g(x_1,\ldots,x_n) = b_m x_1^m + b_{m-1} x_1^{m-1} + \cdots + b_1 x_1 + b_0$, where $a_i, b_i \in \mathbf{R}[x_2, \ldots, x_n]$. The *resultant matrix* $A = A(f, g)$ is defined by the $(k + m) \times (m + k)$ semi-circulant matrix

$$A = \begin{bmatrix} a_0 & 0 & \cdots & b_0 & 0 & \cdots \\ a_1 & a_0 & \cdots & b_1 & b_0 & \cdots \\ \vdots & \vdots & \cdots & \vdots & \vdots & \cdots \\ a_i & a_{i-1} & \cdots & b_i & b_{i-1} & \cdots \\ \vdots & \vdots & \cdots & \vdots & \vdots & \cdots \end{bmatrix}. \qquad (45)$$

The matrix $A$ has $m$ columns of the coefficients of $f$, and $k$ columns of the coefficients of $g$. The *resultant* of a pair of polynomials $f(x_1,..,x_n)$, $g(x_1,..,x_n) \in \mathbf{R}[x_1,..,x_n]$ is defined as the determinant of the corresponding matrix, viz, $\mathrm{Res}(f, g, x_1) = \det(A)$.

The *discriminant* of a monic polynomial is the expression $\mathrm{disc}(f) = (-1)^{k(k-1)/2} \mathrm{Res}(f, f', x_1)$, where $f'$ is the derivative of $f$. The discriminant is widely used to identify polynomials with multiple roots. Similarly, the vanishing resultant $\mathrm{Res}(f, g, x_1) = 0$ implies that the polynomials $f$ and $g$ have a common factor.

There are $n$ (distinct) resultants of a pair of polynomials $f(x_1,..,x_n)$, $g(x_1,..,x_n) \in \mathbf{R}[x_1,..,x_n]$, one for each variable $x_i$, $i = 1, 2, \ldots, n$. The resultant effectively eliminates a variable in the system of polynomial equations $f(x_1,..,x_n) = 0$, $g(x_1,..,x_n) = 0$.

Some useful and elementary properties of the resultants are stated below.

***Theorem* 21.** Let $f(x_1,..,x_n)$, $g(x_1,..,x_n) \in \mathbf{R}[x_1,..,x_n]$ be nonconstant polynomials. Then
(i) The resultant $\mathrm{Res}(f, g, x_1) \in \mathbf{R}[x_2,..,x_n]$.
(ii) The reverse resultant $\mathrm{Res}(g, f, x_1) = (-1)^{km} \mathrm{Res}(f, g, x_1)$.
(iii) Multiplicative property $\mathrm{Res}(g, f, x_1) = \mathrm{Res}(f_1, g, x_1)\mathrm{Res}(f_2, g, x_1)$ if $f = f_1 f_2$.
(iv) $f(x_1,..,x_n)s + g(x_1,..,x_n)t = \mathrm{Res}(f, g, x_1)$, where $s, t \in \mathbf{R}[x_1,..,x_n]$.

A more practical way of computing the resultant of a pair of relatively prime polynomials is by means of the Euclidean algorithm. To realize this rearrange the equation as

$f(x_1,..,x_n)S + g(x_1,..,x_n)T = 1$, where the resultant is the denominator of the fractions

$$S = \frac{s(x_1,\ldots,x_n)}{\mathrm{Res}(f,g)}, \quad T = \frac{t(x_1,\ldots,x_n)}{\mathrm{Res}(f,g)}. \qquad (46)$$

***Example* 22.** The resultant of the system of equations

$$\begin{aligned} a(x, y, z) &= a_3 xz + a_2 x + a_1 z + a_0 = (a_3 x + a_1)z + a_2 x + a_0, \\ b(x, y, z) &= b_3 yz + b_2 y + b_1 z + b_0 = (b_3 y + b_1)z + b_2 y + b_0, \end{aligned} \qquad (47)$$

in three variables effectively eliminates one of the variable, one of the 3 possible resultant is

$$\mathrm{Res}(a,b,z) = \det \begin{bmatrix} a_3 x + a_1 & b_3 x + b_1 \\ a_2 x + a_0 & b_2 x + b_0 \end{bmatrix} = c_3 xy + c_2 x + c_1 y + c_0. \qquad (48)$$

### 4.3 Polynomial Equations

The application of lattice reduction theory to polynomial factorization was pioneered by [L], and the application to the determination of the roots of nonlinear algebraic equations seems to have been the works of [HD] and [VE]. However, the application to linear equations both modular and over the integers has an earlier beginning, see [BJ] and [LR]. Later the technique for nonlinear algebraic equations was significantly improved in [CR]. Specifically, the range of the roots $x$ of the polynomial equation $f(x) = a_d x^d + \cdots + a_1 x + a_0$ modulo $N$ that can be determined in deterministic logarithmic time complexity was extended from $|x| \leq N^{2/(d(d+1))}$ to $|x| \leq N^{1/d-\varepsilon}$, $\varepsilon > 0$. This is accomplished by replacing the original basis of the polynomials lattice.

Since then over a dozen papers and a few dissertations have been published on the applications of lattice reductions theory to integer factorization and cryptography. These more recent works have simplified the theory and its practical aspect.

The height of a polynomial $f(x, y) = \sum_{0 \leq i, j \leq d} a_{i,j} x^i y^j \in \mathbb{Z}[x, y]$ of maximum degree $\deg(f) = d$ in $x$ and $y$ is defined by the supremum norm $\| f(x,y) \| = \| f(x,y) \|_\infty = \max\{ | a_{0,0} |, | a_{0,1} |, \ldots, | a_{d,d} | \}$. Likewise, the height of a rational number $r/s$ is defined by $\| r/s \|_\infty = \max\{ | r |, | s | \}$. The weight of a polynomial is defined by $w(f) = \#\{ | a_{i,j} | \neq 0 : 0 \leq i, j \leq d \}$. The standard norm is defined by the relation $\| f(x, y) \|_2 = \sqrt{a_{0,0}^2 + a_{0,1}^2 + \cdots + a_{d,d}^2}$.

***Theorem 23.*** ([ST]) Let $a(x_1, \ldots, x_n)$ and $b(x_1, \ldots, x_n)$ be two non-zero polynomials over $\mathbb{Z}$ of maximum degree $d$ in each variable separately such that $b(x_1, \ldots, x_n)$ is a multiple of $a(x_1, \ldots, x_n)$ in $\mathbb{Z}[x_1, \ldots, x_n]$. Then $\| b \|_2 \geq 2^{-(d+1)^n + 1} \| a \|_\infty$

***Theorem 24.*** ([NH]) Let $f(x_1, \ldots, x_n) \in \mathbb{Z}[x_1, \ldots, x_n]$ be a polynomial in $n$-variables with $w > 0$ nonzero terms. Suppose that $f(x_1, \ldots, x_n) \equiv 0 \bmod N$, and $\| f(x_1 X_1, \ldots, x_n X_n) \| < N w^{-1/2}$, where $| x_i | < X_i$ for $i = 1, 2, \ldots, n$. Then $f(x_1, \ldots, x_n) = 0$ holds over the integers $\mathbb{Z}$.

***Theorem 25.*** ([CR]) Let $f(x, y) \in \mathbb{Z}[x, y]$ be irreducible with maximum degree $d$ in $x, y$ separately. Let $|X|, |Y|$ be upper bounds on the desired integer solution $(x_0, y_0)$ and let $W = \max\{ | a_{i,j} X^i Y^j | 0 \leq i, j \leq d \}$. If $XY \leq W^{2/3d}$ then all integer pairs $(x_0, y_0)$ such that $f(x_0, y_0) = 0$, $| x_0 | \leq X$ and $| y_0 | \leq Y$ can be found in time polynomial in $\log W$ and $2^d$.

Although Theorem 23 calls for an irreducible polynomial, Theorem 24 seems to imply that reducible polynomials work as well as irreducible polynomials whenever $\| b \|_2 < 2^{-(d+1)^n + 1} \| a \|_\infty$ holds.

A polynomial $f(x, y, z)$ is said to be defined over a subset of monomials $M$ if $f(x, y, z)$ can be written as linear combination of monomials in $M$. Let $S$ be another non-empty set and $g, h$ be two polynomials such that $h(x, y, z) = g(x, y, z) f(x, y, z)$. The ordered pair $(S, M)$ is said to be admissible for *f* if the property "*h* defined over *M*" is equivalent with "*g* defined over *S*", see [BA] for an extended discussion.

The integer $d_1$, $d_2$ and $d_3$ are the maximum degree of the polynomial $f(x, y, z)$ in the variables $x$, $y$, and $z$ respectively. The integers $s_1$, $s_2$ and $s_3$ are defined by the sums

$$s_1 = \sum_{(i,j,k) \in M \setminus S} i \ , \quad s_2 = \sum_{(i,j,k) \in M \setminus S} j \ , \quad \text{and} \quad s_3 = \sum_{(i,j,k) \in M \setminus S} k \ . \tag{49}$$

Let $f(x, y, z)$ be an irreducible polynomial of $\mathbb{Z}[x, y, z]$, and let $(x_0, y_0, z_0)$ be a small root over the integers such that $|x_0| < X$, $|y_0| < Y$ and $|z_0| < Z$.

**Theorem 26.** ([BA]) If $S$ and $M$ are admissible sets for $f(x, y, z)$, then an algebraically independent polynomial $g(x, y, z)$ which has $(x_0, y_0, z_0)$ as a root over the integers can be found in polynomial time, provided that

$$X^{s_1} Y^{s_2} Z^{s_3} < W^s 2^{-(6+c)(d_1^2+d_2^2+d_3^2)s} , \tag{50}$$

where it is assumed that $(m - s)^2 \leq cs(d_1^2 + d_2^2 + d_3^2)$ for some constant $c$.

**Corollary 27.** Let $f(x, y, z) \in \mathbb{Z}[x,y,z]$ be an irreducible polynomial of maximal degree $d > 0$ in $x$, $y$, and $z$, (or of total degree $d > 0$ in $x$, $y$, and $z$) and let $\| f(xX, yY, zZ) \| = W$ be the height of $f$. Suppose that there exists a triple $(x_0, y_0, z_0)$ such that $f(x_0, y_0, z_0) = 0$, where $0 \leq |x_0 y_0 z_0| < W^{2/3d}$, and $0 \leq |z| \leq O((\log N)^B)$, $B > 0$ constant. Then the solution $(x_0, y_0, z_0)$ can be determined in deterministic logarithmic time.

Proof: By Theorem 26, there exists an algebraically independent polynomial $g(x, y, z)$ which has a common small solution $(x_0, y_0, z_0)$ and can be found in deterministic logarithmic time. To verify the ranges of the variables $x$ and $y$, observe that since $O((\log N)^B) = o(N^\varepsilon)$ for any arbitrarily small $\varepsilon > 0$, the contribution of the term $Z^{s_3}$ to the lattice reduction inequality (50) is negligible. Thus for all practical purpose (50) collapses to a system of two variables, and by Theorem 25, it follows that the ranges of the variables $x$ and $y$ satisfy the conditions

$$0 \leq |x_0| \leq X, \ 0 \leq |y_0| \leq Y, \text{ and } XY \leq W^{2/3d} \text{ or } XY \leq W^{1/d},$$

where of $d > 0$ is the maximal degree in $x$, $y$, and $z$, or of total degree $d > 0$ in $x$, $y$, and $z$.

The root $x_0$ or $y_0$ is resolved by an exhaustive search over $0 \leq |z_0| \leq O((\log N)^B)$ using the resultant $\text{Res}_x(f(x, y, z_0), g(x, y, z_0), y)$ or $\text{Res}_y(f(x, y, z_0), g(x, y, z_0), x)$. ∎